\documentclass[journal]{IEEEtran}
\usepackage[]{algorithm2e}
\usepackage[]{algorithmic}

\usepackage[pdftex]{graphicx} 
\usepackage{epstopdf} 
\usepackage{cite}
\usepackage{subfig}
\usepackage{multirow}
\usepackage{amsmath}
\usepackage{amsfonts}
\usepackage{color}
\usepackage{smartdiagram}
\usetikzlibrary{matrix,arrows.meta}
\definecolor{orange}{rgb}{1,0.5,0}
\definecolor{armygreen}{rgb}{0.35, 0.50, 0.30}

\usepackage{layout}
\usepackage{hhline}
\usepackage{siunitx}

\usepackage{nomencl}
\makenomenclature

\usepackage{cases}
\usepackage{dblfloatfix}
\usepackage{tabularx,booktabs}
\newcolumntype{C}{>{\centering\arraybackslash}X} 
\usepackage{lipsum}
\usepackage{float}
\usepackage{hyperref}
\usepackage{tabu}
\usepackage{boldline,multirow}
\usepackage[T1]{fontenc}
\usepackage[11pt]{moresize}
\usepackage{array}
\newcolumntype{P}[1]{>{\centering\arraybackslash}p{#1}}
\usepackage{tikz}
\makeatother
\usepackage{verbatim}
\usepackage{arydshln}
\usetikzlibrary{shapes,arrows, chains}
\usetikzlibrary{positioning, calc}
\usepackage{caption}
\usepackage{bm}
\usepackage{breqn}
\usepackage{tikz,pgfplots}
\usepackage{epsfig}
\usepackage{eurosym}

\usepackage{breqn}
\usetikzlibrary{mindmap,backgrounds}

\usepackage{mathtools}
\DeclarePairedDelimiter\ceil{\lceil}{\rceil}
\DeclarePairedDelimiter\floor{\lfloor}{\rfloor}

\tikzfading[name=priorityarrowfadingdown,
top color=transparent!5,
bottom color=transparent!80
]

\tikzset{priority
arrow fill/.style={
  fill=gray,
  path fading=priorityarrowfadingdown
  }
}

\makeatletter
\NewDocumentCommand{\smartdiagramx}{r[] m}{
    \StrCut{#1}{:}\diagramtype\option
    \IfStrEq{\diagramtype}{priority descriptive diagram}{
        \pgfmathparse{subtract(\sm@core@priorityarrowwidth,\sm@core@priorityarrowheadextend)}
        \pgfmathsetmacro\sm@core@priorityticksize{\pgfmathresult/2}
        \pgfmathsetmacro\arrowtickxshift{(\sm@core@priorityarrowwidth-\sm@core@priorityticksize)/2}
        \begin{tikzpicture}[every node/.style={align=center,let hypenation}]
        \foreach \smitem [count=\xi] in {#2}{\global\let\maxsmitem\xi}
        \foreach \smitem [count=\xi] in {#2}{%
            \edef\col{\@nameuse{color@\xi}}
            \node[description,drop shadow](module\xi)
            at (0,0+\xi*0.65*\sm@core@descriptiveitemsysep) {\smitem};
            \draw[line width=\sm@core@prioritytick,\col]
            ([xshift=-\arrowtickxshift pt]module\xi.base west)--
            ($([xshift=-\arrowtickxshift pt]module\xi.base west)-(\sm@core@priorityticksize pt,0)$);
        }%
        \coordinate (A) at (module1);
        \coordinate (B) at (module\maxsmitem);
        \CalcHeight(A,B){heightmodules}
        \pgfmathadd{\heightmodules}{\sm@core@priorityarrowheightadvance}
        \pgfmathsetmacro{\distancemodules}{\pgfmathresult}
        \pgfmathsetmacro\arrowxshift{\sm@core@priorityarrowwidth/2}
        \begin{pgfonlayer}{background}
        \node[priority arrow,rotate=180,transform shape] at ([xshift=-\arrowxshift pt]module\maxsmitem.north west){};
        \end{pgfonlayer}
        \end{tikzpicture}
    }{}
}%

\makeatother
\begin{document}
\bstctlcite{IEEEexample:BSTcontrol}
\title{Global Optimization of Offshore Wind Farm Collection Systems}
\author{Juan-Andr\'es~P\'erez-R\'ua,
        Mathias Stolpe, Kaushik~Das,~\IEEEmembership{Member,~IEEE}
        and~Nicolaos~A.~Cutululis,~\IEEEmembership{Senior Member,~IEEE}
\thanks{Juan-Andr\'es~P\'erez-R\'ua, Mathias Stolpe, Kaushik Das, and Nicolaos A. Cutululis are with the DTU Wind Energy, Technical University of Denmark, Frederiksborgvej 399, 4000 Roskilde, Denmark (e-mail: juru@dtu.dk)}
}
\markboth{}%
{Shell \MakeLowercase{\textit{et al.}}: Bare Demo of IEEEtran.cls for IEEE Journals}
\maketitle
\begin{abstract}
A mathematical program for global optimization of the cable layout of Offshore Wind Farms (OWFs) is presented. The model consists on a Mixed Integer Linear Program (MILP). Modern branch-and-cut solvers are able to solve large-scale instances, defined by more than hundred Wind Turbines (WTs), and a reasonable number of Offshore Substations (OSSs). In addition to the MILP model to optimize total cable length or initial investment, a pre-processing strategy is proposed in order to incorporate total electrical power losses into the objective function. High fidelity models are adapted to calculate cables current capacities, spatial currents. The MILP model is embedded in an iterative algorithmic framework, solving a sequence of problems with increasing search space size. The search space is defined as a set of underlying candidate arcs. The applicability of the method is illustrated through $10$ case studies of real-world large-scale wind farms. Results show that: (i) feasible points are obtained in seconds, (ii) points with an imposed maximum tolerance near the global optimum are calculated in a reasonable computational time in the order of hours, and (iii) the proposed method compares favorably against state-of-the-art method available in literature.
\end{abstract}
\IEEEpeerreviewmaketitle
\section*{Nomenclature}

\subsection*{Acronyms}
\begin{description}
\item[OWF(s)]   \hspace{0.4cm} Offshore Wind Farm(s).
\item[WT(s)]    \hspace{0.4cm} Wind Turbines(s).
\item[OSS(s)]   \hspace{0.4cm} Offshore Substation(s).
\item[BIP]     \hspace{0.4cm} Binary Integer Programming.
\item[MILP]     \hspace{0.4cm} Mixed Integer Linear Programming.
\item[MIQP]     \hspace{0.4cm} Mixed Integer Quadratic Programming.
\item[MINLP]     \hspace{0.4cm} Mixed Integer Nonlinear  Programming.
\item[SCETM]    \hspace{0.4cm} Single-Core Equivalent Thermal Model.
\item[NP] \hspace{0.4cm} Non-Deterministic Polynomial.
\item[L]    \hspace{0.4cm} Length.
\item[LP]    \hspace{0.4cm} Length plus total Power losses.
\item[I]    \hspace{0.4cm} Initial investment.
\item[IP]    \hspace{0.4cm} Investment plus total Power losses.
\item[WDS]    \hspace{0.4cm} West of Duddon Sands.
\item[TH]    \hspace{0.4cm} Thanet.
\item[LA]    \hspace{0.4cm} London Array.
\item[HR1]    \hspace{0.4cm} Horns Rev 1.
\item[O]    \hspace{0.4cm} Ormonde.
\item[DT]    \hspace{0.4cm} DanTysk.
\end{description}
\subsection*{Parameters (non-sets)}
\begin{description}
\item[$n_w$] \hspace{0.4cm} Number of wind turbines.
\item[$n_o$] \hspace{0.4cm} Number of offshore substations.
\item[$m$] \hspace{0.4cm} Total number of years.
\item[$a_{ij}$] \hspace{0.4cm} Arc connecting point $i$ to $j$ $(i,j)$.
\item[$d_{ij}$] \hspace{0.4cm} Euclidean norm for arc $(i,j)$.
\item[$u_t$] \hspace{0.4cm} Capacity of cable $t$ in number of wind turbines.
\item[$c_{c_{t}}$] \hspace{0.4cm} Metric capital cost of cable $t$.
\item[$c_{p_{t}}$] \hspace{0.4cm} Metric installation cost of cable $t$.
\item[$I_t$] \hspace{0.4cm} Steady-state rated current of cable $t$.
\item[$V_n$] \hspace{0.4cm} Nominal line-to-line voltage of the system.
\item[$P_n$] \hspace{0.4cm} Nominal power of the wind turbines.
\item[$S_{r_{t}}$] \hspace{0.4cm} Nominal power of cable $t$.
\item[$\vec{\gamma}_t$] \hspace{0.4cm} Propagation constant of cable $t$.
\item[${\vec{Z_c}}_{t}$] \hspace{0.4cm} Characteristic impedance of cable $t$.
\item[$\vec{z}_t$] \hspace{0.4cm} Metric series impedance of cable $t$.
\item[$\vec{y}_t$] \hspace{0.4cm} Metric admittance of cable $t$.
\item[$\vec{I}_{ij,t}^{\,k}$] \hspace{0.4cm} Nominal phasor current of arc $(i,j)$ using cable
\item[] \hspace{0.4cm} $t$, when $k$ wind turbines are connected. 
\item[$S_{ij,t}^{k}$] \hspace{0.4cm} Nominal power of arc $(i,j)$ using cable $t$, when
\item[] \hspace{0.4cm} $k$ wind turbines are connected.
\item[$p^{\omega}$] \hspace{0.4cm} Power produced by a wind turbine at $\omega$ hour-slot.
\item[$f^{\omega,k}_{ij}$] \hspace{0.4cm} Power flow in arc $(i,j)$ at $\omega$ hour-slot, 
\item[] \hspace{0.4cm} when $k$ wind turbines are connected. 
\item[$\vec{I}_{ij,t}^{\,\omega,k}$] \hspace{0.4cm} Phasor current through arc $(i,j)$ using cable $t$, at
\item[] \hspace{0.4cm} $\omega$ hour-slot, when $k$ wind turbines are connected. 
\item[$\lambda_{1}$] \hspace{0.4cm} Screen losses factor. 
\item[$\lambda_{2}$] \hspace{0.4cm} Armouring losses factor. 
\item[$W_{d_{t}}$] \hspace{0.4cm} Metric dielectric loss of cable $t$.
\item[$R_{t}$] \hspace{0.4cm} Metric electrical resistance of cable $t$.
\item[$l^{\mu,k}_{ij,t}$] \hspace{0.4cm} Annual total power losses through arc $(i,j)$ using
\item[] \hspace{0.4cm} cable $t$, at year $\mu$, when $k$ wind turbines are co-
\item[] \hspace{0.4cm} nnected. 
\item[$U$] \hspace{0.4cm} Capacity of the biggest cable available given as \item[] \hspace{0.4cm} maximum number of supportable wind turbines.
\item[$r$] \hspace{0.4cm} Discount rate.
\item[$c_{ij}^{k}$] \hspace{0.4cm} Metric cost of arc $(i,j)$, when $k$ wind turbines
\item[] \hspace{0.4cm} are connected.
\item[$\phi$] \hspace{0.4cm} Maximum number of feeders per offshore substa- \item[] \hspace{0.4cm} tion.
\item[$\eta$] \hspace{0.4cm} Loading symmetry coefficient for offshore subs- \item[] \hspace{0.4cm} tations.
\item[$\upsilon$] \hspace{0.4cm} Number of wind turbines arcs set to a wind turbi-
\item[] \hspace{0.4cm} ne. 
\item[$\upsilon_f$] \hspace{0.4cm} Number of wind turbines arcs set to a wind turbi-
\item[] \hspace{0.4cm} ne for the feasibility problem.  
\item[$\upsilon_{f_{min}}$] \hspace{0.4cm} Algorithm parameters for the feasibility problem.
\item[$\upsilon_{f_{\delta}},\upsilon_{f_{max}}$] 
\item[$\upsilon_o$] \hspace{0.4cm} Number of wind turbines arcs set to a wind turbi-
\item[] \hspace{0.4cm} ne for the global optimization problem.  
\item[$\upsilon_{o_{min}}$] \hspace{0.33cm} Algorithm parameters for the global opt. problem.
\item[$\upsilon_{o_{\delta}},\upsilon_{o_{max}}$] 
\item[$\epsilon$] \hspace{0.4cm} Required relative optimality gap.
\end{description}
\subsection*{Parameters (sets)}
\begin{description}
\item[$\boldsymbol{N_{o}}$] \hspace{0.4cm} Set of offshore substations.
\item[$\boldsymbol{N_{w}}$] \hspace{0.4cm} Set of wind turbines.
\item[$\boldsymbol{N}$] \hspace{0.4cm} Set of offshore substations and wind turbines.
\item[$G$] \hspace{0.4cm} Weighted directed graph.
\item[$\boldsymbol{A}$] \hspace{0.4cm} Set of available arcs.
\item[$\boldsymbol{D}$] \hspace{0.4cm} Set of arcs weights.
\item[$\boldsymbol{T}$] \hspace{0.4cm} Set of available cables.
\item[$\boldsymbol{U}$] \hspace{0.4cm} Set of cables capacities in wind turbines number.
\item[$\boldsymbol{C_{c}}$] \hspace{0.4cm} Set of cables capital expenditures costs.
\item[$\boldsymbol{C_{p}}$] \hspace{0.4cm} Set of cables installation costs.
\item[$\boldsymbol{\Omega}^{\mu}$] \hspace{0.4cm} Set of hours-slot for a year $\mu$.
\item[$\boldsymbol{M}$] \hspace{0.4cm} Set of operational years.
\item[$G_{r}$] \hspace{0.4cm} First reduced graph.
\item[$\boldsymbol{A_{r}}$] \hspace{0.4cm} Set of first reduced arcs.
\item[$\boldsymbol{\chi}$] \hspace{0.4cm} Set of crossing pairs arcs.
\item[$\boldsymbol{\Upsilon}_i$] \hspace{0.4cm} Set of wind turbines connected to $i$.
\item[$G^{'}_{r}$] \hspace{0.4cm} Second reduced graph.
\item[$\boldsymbol{A^{'}_{r}}$] \hspace{0.4cm} Set of second reduced arcs.
\end{description}
\subsection*{Variables}
\begin{description}
\item[$x_{ij}$] \hspace{0.4cm} Binary variable to activate arc $(i,j)$.
\item[$x_{ij,t}$] \hspace{0.4cm} Binary variable to select optimum cable type $t$ 
\item[] \hspace{0.4cm} for arc $(i,j)$.
\item[$y^{k}_{ij}$] \hspace{0.4cm} Binary variable to activate arc $(i,j)$, when $k$ wind
\item[] \hspace{0.4cm} turbines are connected.
\item[$\sigma_i$] \hspace{0.4cm} Integer variable of number of wind turbines co-
\item[] \hspace{0.4cm} nnected to offshore substation $i$.
\end{description}
\subsection*{Optimization output (non-sets)}
\begin{description}
\item[$\epsilon_{k_{f}}$] \hspace{0.4cm} Calculated gap at iteration $k_{f}$.
\item[$\epsilon_{k_{o}}$] \hspace{0.4cm} Calculated gap at iteration $k_{o}$.
\item[$\epsilon_{k_g}$] \hspace{0.4cm} Recalculated gap at iteration $k_{g}$.
\end{description}
\subsection*{Optimization output (sets)}
\begin{description}
\item[$\boldsymbol{I}$] \hspace{0.4cm} Set of first feasible solution.
\item[$\boldsymbol{O}_{k_{o}}$] \hspace{0.4cm} Set of feasible solution at iteration $k_{o}$.
\item[$\boldsymbol{\Gamma}_{k_{o}}$] \hspace{0.4cm} Set of candidate arcs at iteration $k_{o}$.
\item[$\boldsymbol{Z}_{k_{o}}$] \hspace{0.4cm} Set of active variables $x_{ij}=1$ of the problem
\item[] \hspace{0.4cm} defined in the iteration $k_o$.
\end{description}

\subsection*{Subscripts}
\begin{description}
\item[$i$] \hspace{0.4cm} Element in the set $\boldsymbol{N}$.
\item[$j$] \hspace{0.4cm} Element in the set $\boldsymbol{N}$.
\item[$ij$] \hspace{0.4cm} Arc $(i,j)$ with tail at $i$ and head at $j$.
\item[$k_{f}$] \hspace{0.4cm} Iteration for the feasibility problem.
\item[$k_{o}$] \hspace{0.4cm} Iteration for the global optimization problem.
\item[$k_{g}$] \hspace{0.4cm} Iteration for the general problem.
\item[$t$] \hspace{0.4cm} Cable type in the set $\boldsymbol{T}$.
\end{description}

\subsection*{Superscripts}
\begin{description}
\item[$k$] \hspace{0.4cm} Number of turbines connected in $(i,j)$.
\item[$\omega$] \hspace{0.4cm} Hour-slot $\in$ $\boldsymbol{\Omega}^{\mu}$.
\item[$\mu$] \hspace{0.4cm} Year $\in$ $\boldsymbol{M}$.
\end{description}

\section{Introduction} 
\label{introduction}
\IEEEPARstart{O}{ffshore} wind energy represents a backbone technology towards the transition to power systems fully based on renewable energy. \\After the invention, and experimentation age during the 1990s, the commercialization and development period is fundamentally focused on turning this technology into not only an environmentally sustainable paradigm, but also financially competitive compared to other classic, and emergent types of energy generation. \\The share of Offshore Wind Farms (OWFs) has increased almost five times in the last seven years \cite{GWEC2017Global2017}, reaching a globally installed power of nearly \SI{19}{GW}. OWF projects are capital intensive, having large values of operating leverage where the required electrical infrastructure costs can raise up to 15\% compared to the total system costs \cite{Sun2012TheDevelopment}.\\ 
The electrical collection system is the set of electrical infrastructure components (AC submarine cables, switchgears, transformers, protection, and control units, etc.). This is required to interconnect the Wind Turbines (WTs) with each other and the Offshore Substation (OSS), guaranteeing an effective, reliable, and efficient collection of energy to the export infrastructure. \\Between 2018 and 2028 more than \SI{19000}{\km} of cables for collection systems are prognosed to be installed in UK only, with an estimated worth of \textsterling5.36bn \cite{RenewableUK2018RenewableUKIntelligence}. Economies of scale pushes the development of large-scale OWFs, having more than 80 WTs while increasing their rated power. \\
The collection system design and optimization problem has been studied with accentuated focus in the last ten years \cite{Perez-Rua2019ElectricalReview,Lumbreras2013OffshoreReview}. Finding the global optimum of this problem is generally NP-hard \cite{Jothi2005ApproximationDesign}. Four big clusters of methods for tackling this problem can be established: heuristics, metaheuristics, global optimization with mathematical formulations, and hybrids, such as matheuristics.\\
Global optimization encompasses a large set of different alternatives to model the cable layout problem, like Binary Integer Programming (BIP) \cite{Cerveira2016OptimalCase},  Mixed Integer Linear Programming (MILP) \cite{Bauer2015TheProblem, Fischetti2018OptimizingLosses,Fischetti2018MixedRouting,Berzan2011AlgorithmsFarms,Wedzik2016AOptimization}, MILP with decomposition techniques for stochastic programming \cite{Lumbreras2013AProblems,Lumbreras2013OptimalStrategies}, Mixed Integer Quadratic Programming (MIQP) \cite{Hertz2017DesignAvailable,Banzo2011StochasticFarms}, and Mixed Integer Non-Linear Programming (MINLP) \cite{Chen2016CollectorFarms,Pillai2015OffshoreOptimization}. \\Important advances on mathematical modelling are provided in \cite{Cerveira2016OptimalCase}, but disregarding fundamental practical considerations for OWFs, such as cables crossings, multiple OSSs, maximum number of feeders per OSS, wind power stochasticity, among others. The application focuses on global optimization for medium-scale OWFs (with the largest problem instance of $57$ units). Similarly, \cite{Berzan2011AlgorithmsFarms} and \cite{Wedzik2016AOptimization} are tackling small-scale projects (with less than or equal to $30$ WTs), without proposing any strategies for scaling the applications to large-scale problems and, as in  \cite{Cerveira2016OptimalCase}, do not include practical considerations in their modelling. \\Larger instances are designed to optimality in \cite{Bauer2015TheProblem}, however, following a Planar Open Vehicle Routing Problem approach (no branching), restricting the cables set size, and ignoring power losses. Large-scale OWFs are tackled in \cite{Fischetti2018OptimizingLosses}  and \cite{Fischetti2018MixedRouting}, combining a MILP flow-based model with up to four heuristics, considering power losses, and other practical applications. 
\\The works \cite{Lumbreras2013OptimalStrategies,Lumbreras2013AProblems}, and \cite{Banzo2011StochasticFarms} provide remarkable advances on stochastic optimization for problems in this context. Different stochastic scenarios are supported, accounting for wind power variability, and cables failure. Distinctive theoretical strategies to accelerate convergence are applied and compared. Nonetheless, some simplifications are incorporated, such as iterative processes to calculate total electrical power losses, combined with experts analysis to estimate their impact on the cable layout. Case studies are limited to small-scale applications.\\By means of explicit formulation  of electrical losses economic costs in the model, a MINLP program is proposed in \cite{Pillai2015OffshoreOptimization}, where, additionally, clustering algorithms are used prior the execution of the program into a commercial solver. Finally, losses can be also computed in the accurate quadratic form as in \cite{Hertz2017DesignAvailable}, but at the expense of a decrease in computational efficiency.
\\Each of these mathematical formulations impose certain limitations about the physics modelling options. For instance, using flow-based MILP makes it more difficult to include the quadratic active power losses explicitly into the objective function. The commonly used power flow equations solved with e.g. the Newton-Raphson method cannot be considered in MILP or MIQP formulations. \\The side effects of more flexible modelling formulations are the compromise of solver functionality and performance. The proper balance between solution method and complexity on modelling, represents one of the main challenges for the OWF developer, and trade-offs have to be adopted within certain assumptions. As a generalization, linear-based formulations are computationally more efficient than quadratic or non-linear.
\\Since the trend in OWFs is to deploy large-scale projects, focus is directed into this aspect. To the best of the authors' knowledge, only \cite{Bauer2015TheProblem}, \cite{Fischetti2018OptimizingLosses}, \cite{Pillai2015OffshoreOptimization}, and \cite{Klein2017Obstacle-awareLayouts} have tackled collection systems for large-scale OWFs with $80$ - $100$ WTs using global optimization. Only one cable type and no-branching at WTs nodes is considered by Bauer et. al. \cite{Bauer2015TheProblem}, where heuristics have also been proposed. A matheuristic framework is developed by Fischetti et. al. \cite{Fischetti2018OptimizingLosses}, including the total electrical power losses in the objective function using pre-processing strategies, and flow formulation. A MINLP flow formulation is used in \cite{Pillai2015OffshoreOptimization}, having sets of linear constraints but with explicit linearized losses inclusion in the objective; a fixed neighboring search is implemented, and pre-clustering of WTs to OSS is required, limiting artificially the search space, leading to sub-optimal points. No losses have been included in the optimization in \cite{Klein2017Obstacle-awareLayouts}, and for large-scale problems, computational experiments point out that feasible solutions may not be found. Additionally, only a single cable is supported.\\
The main contributions of this work are: (i) development, testing, and application of a mathematical model to quickly find feasible points for large-scale OWF instances (with more than $100$ WTs). This is possible by proposing a MILP model with reduced constraints and variables than previously used flow formulations for collection systems in OWFs. (ii) combination of an algorithmic framework with mathematical formulation for obtaining global optimum solution points (or near to it) in reasonable computational time. Additionally, improvements on the complexity and fidelity of modelling aspects for integration with a linear formulation, such as total power losses calculated using time series, and capacitive currents, have been considered. The method does not required pre-clustering of WTs to OSSs, as this is tackled intrinsically in the model.\\In Section \ref{modelling}, the modelling aspects are explained in details; followed by Section \ref{model}, where the optimization model is formulated including objective function and constraints definition. The whole framework description is presented in Section \ref{framework}, where the full mathematical formulation is compacted into a single main iterative algorithm. Computational experiments are performed in Section \ref{experiments}, and the work is finalized with the conclusions in Section \ref{conclusion}.
\section{Modelling aspects}
\label{modelling}
\subsection{Problem graph representation}
\label{graph}
The aim of the optimization is to design the cable layout of the collection system for OWFs, i.e., to interconnect the $n_w$ WTs to the available OSSs, $n_o$. \\In graph theory, a OWF is represented by a forest with roots at the OSSs. A few OWF developers opt for systems with loops. The inclusion of loops may increase reliability, but raises other technical challenges as cables sizing, and generally these systems are designed following heuristic and experts rules. Radial layouts are the common practice in the industry, since array cables failures are rather rare as they are buried in the seabed \cite{Warnock2019FailureSystems}. Therefore, the desired result is $n_o$ spanning-trees, which minimize the required objective function, while satisfying the operational and topological constraints.\\
Let the OSSs define the set $\boldsymbol{N_{o}} = \left\lbrace1,\cdots,n_o\right\rbrace$. Likewise, for the WTs, let $\boldsymbol{N_{w}} = \left\lbrace n_o + 1 ,\cdots,n_o+n_w\right\rbrace$. In this way, each one of the OSSs and WTs (modelled as points in the space) have associated a unique  identifier $i$, such as $i \in \boldsymbol{N} = \boldsymbol{N_{o}} \cup \boldsymbol{N_{w}}$. The Euclidean norm between the positions of the points $i$ and $j$, is defined as $d_{ij}$. The aforementioned inputs are condensed as a weighted directed graph $G(\boldsymbol{N},\boldsymbol{A},\boldsymbol{D})$, where $\boldsymbol{N}$ represents the vertex set, $\boldsymbol{A}$ the set of available arcs arranged as a pair-set, and $\boldsymbol{D}$ the set of associated weights for each element $a_{ij} \in \boldsymbol{A}$, where $i \in \boldsymbol{N}$ $\wedge$ $j \in \boldsymbol{N}$. For instance, for $a_{ij} = (i,j)$, $d = d_{ij}$, where $d \in \boldsymbol{D}$. In general, $G(\boldsymbol{N},\boldsymbol{A},\boldsymbol{D})$ is a complete directed graph.\\
Additionally, a predefined list of available cables types is required to interconnect the WTs towards the OSSs. Let the set of cables be $\boldsymbol{T}$ and let the capacity of a cable $t \in \boldsymbol{T}$ be $u_t$ measured in terms of number of supportable WTs connected downstream. Hence, let $\boldsymbol{U}$ be the set of capacities sorted as in $\boldsymbol{T}$ (see the definition in Section \ref{capacity}). \\Furthermore, each cable type $t$ has a cost per unit of length, $c_{c_{t}}$, in such a way that $u_t$ and $c_{c_{t}}$ describe a positive correlation, following an exponential regression model. The set of metric capital expenditures is defined as $\boldsymbol{C_{c}}$. Similarly the set of metric installation costs is defined by $\boldsymbol{C_{p}}$.\\
After defining the graph representation of the problem, the underlying variables associated to the calculation of the desired output, are unequivocally established. Let $x_{ij}$ represent a binary variable that is one if the arc between the vertex $i$ and $j$ is selected in the solution, and zero otherwise. Likewise, the binary variable $y^{k}_{ij}$ models the $k$ number of WTs connected downstream from $j$, including the wind turbine at node $j$ (under the condition that $x_{ij} = 1$). Finally, the integer variable $\sigma_i$ represents the number of WTs connected to the OSS $i$.
\subsection{Cable capacity}
\label{capacity}
The current capacity $I_t$ of a cable $t$ is calculated using the model given in \cite{IEC2014IEC-60287-1:Rating}. This method represents the most common industrial practice, and it is based on a Single-Core Equivalent Thermal Model (SCETM) as generalized in \cite{Anders1999CalculationFillers} for single-core and three-core cables.
\begin{equation} \label{eqn:capacity}
u_t =  \floor*{\frac{S_{r_{t}} = \sqrt{3} \cdot V_n \cdot I_t}{P_n}} \quad \forall t \in \boldsymbol{T}
\end{equation}
\begin{equation} \label{eqn:capvector}
\boldsymbol{U} = \left\lbrace 1,\cdots,u_{|\boldsymbol{T}|}\right\rbrace
\end{equation}
The set of cable capacities in terms of number of supportable WTs is defined in (\ref{eqn:capacity}) and (\ref{eqn:capvector}), where $P_n$ represents the nominal power of an individual WT, and $V_n$ the line-to-line nominal voltage of the system. \\The cables capacity must be assessed in terms of calculated maximum conductor temperature, and the maximum allowable value by design (usually considered \SI{90}{\celsius}). In this manuscript, the conservative approach described in \cite{IEC2014IEC-60287-1:Rating} is followed. This model neglects the slow thermal time constant of submarine cables, and the high variability of offshore wind power, since it assumes nominal steady state conditions. On the other hand, this practical assumption brings along as a result a robust design, while avoiding power curtailment at any instant.
\subsection{Arcs nominal power}
\label{nominal}
 The power flow model of a transmission line solving the resultant differential equations considering parameters distributed uniformly throughout the length of the cable \cite{Grainger1994PowerAnalysis}, is implemented with the expressions:
 \begin{equation} \label{eqn:current}
\vec{I}_{ij,t}^{\,k} = \frac{k \cdot P_{n}}{\sqrt{3} \cdot V_{n}} \cdot \cosh{(\vec{\gamma}_t} \cdot d_{ij}) -  \frac{V_{n}}{{\vec{Z_c}}_{t}} \cdot \sinh{(\vec{\gamma}_t} \cdot d_{ij}) 
\end{equation}
\begin{equation} \label{eqn:powerarc}
S_{ij,t}^{k} = \sqrt{3} \cdot V_{n} \cdot |\vec{I}_{ij,t}^{\,k}| 
\end{equation}
 The  model (\ref{eqn:current}) provides the highest complexity on physics modelling possible to apply given the flexibility of the proposed mathematical formulation (see Section \ref{model}). This provides more accuracy for securely determining the cable type in an active arc $x_{ij}$. The characteristic impedance is calculated as ${\vec{Z_c}}_{t} = \sqrt{\vec{z}_t}/\vec{y}_t$, and the propagation constant, $\vec{\gamma}_t = \sqrt{\vec{z}_t} \cdot \vec{y}_t$. The series impedance is represented by $\vec{z}_t$, and the admittance by $\vec{y}_t$. The maximum power flowing through the arc $(i,j)$, given $k$ turbines connected downstream ($y^k_{ij}$), is calculated as per (\ref{eqn:powerarc}), accounting for the worst-case scenario, as the current in the arc is strictly increasing with length, and the value $\vec{I}_{ij,t}^{\,k}$ is calculated at the extreme of it when using cable $t$.

\subsection{Power flow and total power losses}
\label{plosses}
A transportation model, accurate enough for radial systems, is implemented following \cite{Perez-Rua2019ElectricalReview}. As explained in Section \ref{capacity}, the set $\boldsymbol{U}$ assumes nominal steady state conditions, describing a practical conservative criterion.
\begin{gather} \label{eqn:powerflow}
\sum\limits_{i \in \boldsymbol{N}}\sum\limits_{k = 1}^{f(i)} k \cdot y^{k}_{ij} \cdot p^{\omega} - \sum\limits_{i \in \boldsymbol{N_{w}}}\sum\limits_{k = 1}^{f(i)} k \cdot y^{k}_{ji} \cdot p^{\omega} = p^{\omega}  \quad \\  \forall j \in \boldsymbol{N_{w}} \wedge \omega \in \boldsymbol{\Omega}^{\mu} \wedge \mu \in {\boldsymbol{M}} \nonumber
\end{gather}
The transportation model is generalized considering the temporal dimension (\ref{eqn:powerflow}), where $\boldsymbol{M} = \left\lbrace1,\cdots,m \right\rbrace$, is the set of operational years with upper limit in the project lifetime, $\boldsymbol{\Omega}^{\mu}$ is the set of hours-slot for a year $\mu$, and $\omega$ a specific hour-slot in $\boldsymbol{\Omega}^{\mu}$. Let $p^w$ be the power in $\si{MW}$ produced by one WT in that hour-slot, simulated with \cite{Srensen2007PowerFarms}. \\
In this way, let define the auxiliary variable $f^{\omega,k}_{ij} = k \cdot y^{k}_{ij} \cdot p^{\omega}$ as the power flow ($\si{MW}$) in arc $(i,j)$ when $k$ WTs are connected downstream (including the one in $j$), in time instant $\omega$. This ignores wake losses, a reasonable assumption since micrositing optimization techniques are considered already applied \cite{Herbert-Acero2014AFarms}.
\begin{equation} \label{eqn:curren2}
\vec{I}_{ij,t}^{\,\omega,k} = \frac{f^{\omega,k}_{ij} }{\sqrt{3} \cdot V_{n}} \cdot \cosh{(\vec{\gamma}_t} \cdot d_{ij}) -  \frac{V_{n}}{{\vec{Z_c}}_{t}} \cdot \sinh{(\vec{\gamma}_t} \cdot d_{ij}) 
\end{equation}
\begin{gather} \label{eqn:losses}
     l^{\mu,k}_{ij,t} \approx 3 \cdot (1+\lambda_{1}+\lambda_{2}) \cdot \sum\limits_{\omega \in \boldsymbol{\Omega}^{\mu}} \omega \cdot R_{t} \cdot d_{ij} \cdot {|\vec{I}_{ij,t}^{\,\omega,k}}|^{2} \\+ 3 \cdot |\boldsymbol{\Omega}^{\mu}| \cdot \omega \cdot W_{d_{t}} \cdot d_{ij} \nonumber
\end{gather}
Including the capacitive currents, (\ref{eqn:curren2}) expresses the current at the end of the arc (with respect to $i$), with magnitude $|\vec{I}_{ij,t}^{\,\omega,k}|$. The annual total power losses $l^{\mu,k}_{ij,t}$ are calculated with ($\ref{eqn:losses}$). The factor $(1+\lambda_{1}+\lambda_{2})$ accounts for the screen and armouring losses and $W_{d_{t}}$ is the dielectric loss per unit length for the insulation surrounding the conductor in \si{W/m}\cite{IEC2014IEC-60287-1:Rating}, while the constant 3 is for the three-phase system. This value must be scaled in $\si{MWh}$.
\section{Optimization Model}
\label{model}
The proposed optimization model described in this Section is able to cope with an arbitrary number of WTs, $n_w$, and similarly  any reasonable number of OSSs, $n_o$.\\ The underlying mathematical formulation is inspired by the formulations and analysis proposed in \cite{Cerveira2016OptimalCase} and \cite{Gouveia2005TheFormulations} with additional constraints stemming from the nature of the problem, and with the objective to improve its tractability.\\
In Section \ref{graph} the binary variables $x_{ij}$ and $y^{k}_{ij}$, and the integer variable $\sigma_i$ are defined. They refer to an active arc, the number $k$ of WTs connect to that active arc, and the number of WTs connected to OSS $i$, respectively. \\To increase the computational efficiency, the number of variables is reduced as follows. The capacity of the biggest cable is calculated as $U = \max{\boldsymbol{U}}$, therefore the possible maximum value of $k$ for $i \in \boldsymbol{N_{o}}$ is equal to $f(i) = U$, while for $i \in \boldsymbol{N_{w}}$ is $f(i) = U-1$. This acknowledges that the biggest cable available could be only used at maximum capacity when is connected from a OSS.\\
Analogously, the set of variables $x_{ij}$, where $i \in \boldsymbol{N_{w}}$, and $j \in \boldsymbol{N_{o}}$ are intrinsically discarded, considering the nature of the power flow, i.e., the OSSs collects the energy from the WTs and not the other way around. Lastly, since the export system is outside the scope of this article, all the arcs between OSSs are disregarded, i.e., $x_{ij} = 0$ $\forall i \in \boldsymbol{N_{o}}$ $\wedge$ $j \in \boldsymbol{N_{o}}$. \\The graph $G(\boldsymbol{N},\boldsymbol{A},\boldsymbol{D})$ is reduced to $G_{r}(\boldsymbol{N},\boldsymbol{A_{r}},\boldsymbol{D_{r}})$ after this stage. While including the OSSs in the optimization could offer some flexibility, their locations are typically decided already in the development process. Furthermore, their locations are strongly driven by the distance to the onshore connection points, hence is deemed as plausible assumption to consider fixed location of the OSSs.
\subsection{Cost coefficients}
\label{costcoefficients}
Note that the previously defined decision variables $x_{ij}$ and $y^{k}_{ij}$, do not include any information related to the cable type selected in a given arc. \\This is because the cable type selection process is handled in a pre-processing stage, given that all the required data is present, and the task is totally independent to any other part of the desired tree(s) \cite{Cerveira2016OptimalCase}. \\A pre-processing strategy allows integrating more complex power flow and total electrical power losses models, increasing the accuracy without compromising the computational efficiency. Along these lines, the distributed model explained in Section \ref{nominal} is harmonized with the main optimization problem. \\At the same time, the aforementioned point allows for a power flow estimation in a conservative fashion, i.e., overestimating the incoming power flow by neglecting the total power losses downstream.\\ In the pre-processing state, for the case of $y^{k}_{ij}$, the length of the arc is known ($d_{ij}$), and the number of WTs connected by it is also defined (by $k$). No more inputs are required for this task. \\Hence, for each $y^{k}_{ij}$, the sub-problem defined by (\ref{eqn:cableselection}) to (\ref{eqn:cableselection4}) is solved beforehand and independently by enumeration.
\begin{eqnarray}
c^{k}_{ij} = \text{min} & \sum\limits_{t \in \boldsymbol{T}} x_{ij,t} \cdot {\left((c_{c_{t}}+c_{p_{t}}) \cdot d_{ij} + \sum\limits_{	\mu=1}^{m} \frac{l^{\mu,k}_{ij,t} \cdot {c_{e}}}{(1+r)^{\mu}}\right)} \label{eqn:cableselection}\\
\text{s.t.} & \sum\limits_{t \in {\boldsymbol{T}}} x_{ij,t} = 1 \label{eqn:cableselection2}\\
    & x_{ij,t} \cdot (S_{ij,t}^{k} - S_{r_{t}}) \leq 0 \quad  \forall t \in \boldsymbol{T} \label{eqn:cableselection3}\\
    & x_{ij,t} \in \left\lbrace0,1\right\rbrace \quad  \forall t \in \boldsymbol{T} \label{eqn:cableselection4}
\end{eqnarray}
In the objective function (\ref{eqn:cableselection}), the first term ($(c_{c_{t}} + c_{p_{t}}) \cdot d_{ij}$) is for capital expenditure ($c_{c_{t}}$) and installation costs ($c_{p_{t}}$) of cable $t$. \\The values for $c_{c_{t}} \in \boldsymbol{C_{c}}$ are obtained from the exponential regression function, given as
\begin{equation} \label{eqn:cost}
c_{c_{t}} = a_{p_{t}}+b_{p_{t}} \cdot \operatorname{e}^{\left(\frac{c_{p_{t}}S_{r_{t}}}{10^8}\right)^{2}}
\end{equation}
where $a_{p_{t}}$, $b_{p_{t}}$, and $c_{p_{t}}$ are coefficients dependent on the nominal voltage of cable type $t \in {\boldsymbol{T}}$, $S_{r_{t}}$ is the rated power of $t$ in \si{VA} (also depending on the rated line to line voltage level, $V_{n}$, see the definition in Section \ref{capacity}). The cables capital expenses cost function (\ref{eqn:cost}) is extracted from \cite{Lundberg2003ConfigurationParks}, being based on a comprehensive semi-empirical cost survey, applying data fitting techniques.\\
The second term in the summation part of (\ref{eqn:cableselection}) accounts for the discounted cash flow of the economic losses caused by the energy dissipation in the cables, the parameters $m$, $l^{\mu,k}_{ij,t}$ (see Section \ref{plosses}), ${c}_{e}$, and $r$, represents the project lifetime (years), total power losses at year $\mu$ for cable $t$ (\si{MWh}), cost of energy (\EUR{}/\si{MWh}), and discount rate (p.u.) respectively. \\The objective function can be simplified by zeroing any of its terms, such as, if only the total length is minimized (L), the first term in (\ref{eqn:cableselection}) has to be replaced uniquely by $d_{ij}$, while the other term in the same equation is dropped. Similarly, if only the total initial investment (I) is targeted, only the first is kept.\\Therefore, the set of single objectives available in the model are: length (L), length plus total power losses (LP)—monetizing lengths by assuming the same capital and installation costs for all cables types, initial investment (I), and initial investment plus total power losses (IP).\\Likewise, (\ref{eqn:cableselection2}) ensures that exactly one cable type is selected, while (\ref{eqn:cableselection3}) guarantees that the capacity of cable $t$ is not exceeded; $S_{ij,t}^{k}$ is the power through arc $(i,j)$ when $k$ turbines are connected in $j$ using cable $t$ defined in Section \ref{nominal}. The binary variable for selecting a cable type $t$ for arc $(i,j)$ is defined in (\ref{eqn:cableselection4}).\\The sub-problem from (\ref{eqn:cableselection}) to (\ref{eqn:cableselection4}) seeks to find the cable type $t$ to be used for the arc $(i,j)$, which minimizes the objective (\ref{eqn:cableselection}).
\subsection{Objective function}
\label{objective}
After solving the multiple sub-problems related to cable selection and cost evaluation (maximum $U \cdot |\boldsymbol{N}|^{2}$ problems) from (\ref{eqn:cableselection}), a cost value $c_{ij}^{k}$ is associated to each $y_{ij}^{k}$ variable. The linear objective function of the main mathematical model is then 
\begin{equation} \label{eqn:objective}
\min{\sum\limits_{i \in \boldsymbol{N}}\sum\limits_{j \in \boldsymbol{N_{w}}}\sum\limits_{k = 1}^{f(i)}  c_{ij}^{k} \cdot y^{k}_{ij}}
\end{equation}
\subsection{Constraints}
\label{constraints}
In order to present the solution connecting all WTs between each other and to the OSSs, the following constraint is added
\begin{equation} \label{eqn:c0}
\sum\limits_{i \in \boldsymbol{N_{o}}} \sigma_{i} = n_w 
\end{equation}
Constraint (\ref{eqn:c0}) models the full OWF to be divided into multiple disconnected trees (forest) with $\sigma_i$ being the number of WTs associated to a OSS $i$. Hence, the total amount of WTs ($n_w$) are integrated into the electrical system.
\\To guarantee full connectivity in OSS $i$, the next constraint is added
\begin{equation} \label{eqn:c1}
\sum\limits_{j \in \boldsymbol{N_{w}}}\sum\limits_{k = 1}^{f(j)}  k \cdot y^{k}_{ij} = \sigma_{i} \quad  \forall i \in \boldsymbol{N_{o}}
\end{equation}
Note that (\ref{eqn:c0}) and (\ref{eqn:c1}) are combined in the case of only one OSS.\\
To limit the maximum number of feeders per OSS ($\phi$), it is used:
\begin{equation} \label{eqn:c12}
\sum\limits_{j \in \boldsymbol{N_{w}}}\sum\limits_{k = 1}^{f(j)} y^{k}_{ij} \leq \phi \quad  \forall i \in \boldsymbol{N_{o}}
\end{equation}
To simultaneously ensure a tree topology, ensure that only one cable type used per arc, and to define the head-tail convention, the next expression is included into the model
\begin{equation} \label{eqn:c2}
\sum\limits_{i \in \boldsymbol{N}}\sum\limits_{k = 1}^{f(i)} y^{k}_{ij} = 1 \quad  \forall j \in \boldsymbol{N_{w}}
\end{equation}
The flow conservation, which also avoids disconnected solutions, is considered by means of one linear equality per wind turbine
\begin{equation} \label{eqn:c3}
\sum\limits_{i \in \boldsymbol{N}}\sum\limits_{k = 1}^{f(i)} k \cdot y^{k}_{ij} - \sum\limits_{i \in \boldsymbol{N_{w}}}\sum\limits_{k = 1}^{f(i)} k \cdot y^{k}_{ji} = 1 \quad  \forall j \in \boldsymbol{N_{w}}
\end{equation}
The set $\boldsymbol{\chi}$ stores pairs of arcs $\left\lbrace (i,j),(u,v)\right\rbrace$, which are crossing each other. Excluding crossing arcs in the solution is ensured by the simultaneous application of the following linear inequalities
\begin{equation} \label{eqn:c5}
x_{ij} + x_{ji} + x_{uv} + x_{vu} \leq 1 \quad  \forall \left\lbrace (i,j),(u,v)\right\rbrace \in \boldsymbol{\chi}
\end{equation}
 \begin{equation} \label{eqn:c4}
\sum\limits_{k = 1}^{f(i)} y^{k}_{ij} - x_{ij} \leq 0 \quad  \forall (i,j) \in \boldsymbol{A_{r}}
\end{equation}
The no-crossing cables restriction is a practical requirement in order to avoid hot-spots, and potential single-points of failure caused by overlapping cables \cite{Bauer2015TheProblem}. Constraint (\ref{eqn:c5}) exhaustively lists all combinations of crossings arcs, including also the corresponding inverse elements. The constraints in (\ref{eqn:c4}) ensure that no active arcs are crossing or overlapping between each other. These constraints thus link the variables $y^{k}_{ij}$ and $x_{ij}$.\\Cables crossings are detected based on a procedure of slopes evaluation. Two arcs are crossing if the crossing point is inside of the lines, but not if this point is located at the extremes of the lines or beyond in the lines' projections.
\begin{gather} \label{eqn:c6}
- \sum\limits_{i \in \boldsymbol{N}}\sum\limits_{k = v + 1}^{f(i)} \floor*{\frac{k - 1}{v}} \cdot y^{k}_{ij} + \sum\limits_{i \in \boldsymbol{N_{w}}}\sum\limits_{k = v}^{f(i)} y^{k}_{ji} \leq 0 \\ \quad \forall v = \left\lbrace 2,\cdots,U-1\right\rbrace \wedge j \in \boldsymbol{N_{w}} \nonumber
\end{gather}
Constraint (\ref{eqn:c6}) represents a set of valid inequalities, initially proposed in \cite{Cerveira2016OptimalCase}, to tighten the mathematical model. Given an active arc $y^{k}_{ij}$, the maximum number of active arcs rooted in $j$ and connecting $v$ WTs, is expressed by $\floor*{\frac{k - 1}{v}}$, hence the constraint restricts the maximum number feasible arcs, reducing the search space without excluding valid solutions to the problem.
\begin{gather} \label{eqn:c7}
x_{ij} \in \left\lbrace0,1\right\rbrace \quad  y^k_{ij} \in \left\lbrace0,1\right\rbrace \\  \forall (i,j) \in \boldsymbol{A_{r}} \wedge k \in \left\lbrace1,\cdots,f(i)\right\rbrace \nonumber
\end{gather}
\begin{equation} \label{eqn:c8}
 0\leq \sigma_{i}\leq \eta \cdot \ceil*{\frac{n_w}{n_o}} \quad \sigma_{i} \in \mathbb{Z_+} \quad \forall i \in \boldsymbol{N_{o}}
\end{equation}
Constraints (\ref{eqn:c7}) and (\ref{eqn:c8}) define the nature of the formulation by the variables definition, a MILP. \\Note that variables $\sigma_i$ are limited in their upper bounds to avoid uneven loading of OSSs (in case $\eta=1$, otherwise $1<\eta\leq n_o$). Equally rated OSSs bring benefits like design standardization, and decreasing of the dependency upon a single transformation unit for transporting the generated power.\\
To summarize, the complete formulation of the  main MILP model consists of the objective function
(\ref{eqn:objective}) and the constraints defined in (\ref{eqn:c0}) - (\ref{eqn:c8}). \\The base formulation presented so far has a maximum number of binary variables equal to $|\boldsymbol{N}|^{2} + U \cdot |\boldsymbol{N}|^{2}$, integer variables number equal to $|\boldsymbol{N_{o}}|$ (linear in function of $n_o$), and constraints (excluding the crossing constraints and valid inequalities) of $1 + 2 \cdot |\boldsymbol{N_{o}}| + 2 \cdot |\boldsymbol{N_{w}}|$. Flow formulations, such as the one proposed in \cite{Fischetti2018OptimizingLosses}, have more variables ($2 \cdot |\boldsymbol{N}|^{2} + U \cdot |\boldsymbol{N}|^{2}$) and constraints ($|\boldsymbol{N}|^2 + 2 \cdot |\boldsymbol{N_{w}}| + |\boldsymbol{N_{o}}|$); integer and binary variables are quadratic in function of the problem size. This fact along with the addition of valid inequalities may explain why the model from (\ref{eqn:objective}) to (\ref{eqn:c8}) is often more efficient to solve. \\Further simplifications to the model are presented in Section \ref{candidate}.
\section{Optimization framework}
\label{framework}
\subsection{Candidate arcs}
\label{candidate}
The reduced graph $G_{r}$ is obtained after performing the described considerations in the introduction of Section \ref{model}. \\However, given the NP-Hard nature of this problem, which is similar to a Capacitated Minimum Spanning Tree with additional constraints \cite{Gonzalez-Longatt2012OptimalApproach,Berzan2011AlgorithmsFarms}, more reductions are required.\\The limitations for successfully finding feasible points and high quality solutions, using solely mathematical models and commercial solvers, is demonstrated in \cite{Fischetti2018OptimizingLosses}.\\For large-scale OWFs (with more than $100$ WTs) the computing time for robust global optimization solvers generally becomes notoriously long. Likewise, in general, solution times become unpredictable, while very large memory requirements are demanded to build the branch-and-cut tree. Besides, the constraints generation must be done with special care (the full set of crossing constraints has a combinatorial nature) to increase computational efficiency. \\To make the formulation more flexible and implementable, a further operation to the graph $G_{r}$ is proposed. The function $f(i,G_{r},\upsilon)$ calculates the set $\boldsymbol{\Upsilon}_i$, defined as the $\upsilon$-closest WTs to $i$. In other words, it is intuitively considered that a WT will be connected to one of the WTs in its vicinity. Therefore, by systematically applying $f(i,G_{r},\upsilon)$ to each $i \in \boldsymbol{N_{w}}$, the reduced graph $G^{'}_{r}$ is found. The set $\boldsymbol{A^{'}_{r}}$ contains the candidate arcs to the solution of the problem. \\With this strategy, the maximum number of variables is reduced to $|\boldsymbol{N_o}| + (U+1) \cdot |\boldsymbol{N_o}| \cdot |\boldsymbol{N_w}| + U \cdot \upsilon \cdot |\boldsymbol{N_w}|$. Additionally, the number of crossing constraints decreases dramatically as well. \\Overall, the arcs set transformation follows $\boldsymbol{A}\rightarrow\boldsymbol{A_r}\rightarrow\boldsymbol{A^{'}_{r}}$.\\ All $(ij)$ indexed variables and constraints, presented in the Section \ref{model}, must be adapted to this reduction process.

\subsection{The Algorithm}
\label{thealgorithm}
The main algorithm defining the full framework with the mathematical model is presented in Algorithm \ref{alg:al2}.\\
\begin{algorithm}[h]
\small
\caption{The main algorithm}
\label{alg:al2}
\algsetup{linenosize=\normalsize, linenodelimiter=., indent=1em}
\rule[-2ex]{0.95\linewidth}{0.5pt} 
\begin{algorithmic}[1]
\STATE{\textit{$k_{f}\leftarrow1$}}
\FOR{\textit{($\upsilon_{f}=\upsilon_{f_{min}}:\upsilon_{f_{\delta}}:\upsilon_{f_{max}})$}}
\STATE{\textit{$G^{'}_{r}\leftarrow f(i,G_{r},\upsilon_{f})\quad \forall i \in \boldsymbol{N_{w}}$}}
\STATE{\textit{$c_{ij}^{k}=0 \quad \forall (i,j) \in \boldsymbol{A^{'}_{r}} \wedge k \in \left\lbrace1,\cdots,f(i)\right\rbrace$}}
\STATE{\textit{Formulate and solve MILP model from (\ref{eqn:objective}) to (\ref{eqn:c8})}}
    \IF{\textit{(problem is feasible)}}
        \STATE{\textit{Save initial feasible point found: $\boldsymbol{I} = x_{ij} \cup y^k_{ij} \quad \forall (i,j) \in \boldsymbol{A^{'}_{r}} \wedge k \in \left\lbrace1,\cdots,f(i)\right\rbrace$}}
        \STATE{\textit{Break}}
    \ELSE
        \STATE{\textit{$k_{f}\leftarrow k_{f}+1$}}
    \ENDIF
\ENDFOR
\STATE{\textit{$k_{o}\leftarrow1$}}
\FOR{\textit{($\upsilon_o=\upsilon_{o_{min}}:\upsilon_{o_{\delta}}:\upsilon_{o_{max}})$}}
\STATE{\textit{$G^{'}_{r}\leftarrow f(i,G_{r},\upsilon_{o})\quad \forall i \in \boldsymbol{N_{w}}$}}
\STATE{\textit{$\boldsymbol{\Gamma}_{k_{o}} = \boldsymbol{A^{'}_{r}}$}}
\STATE{\textit{Get $c_{ij}^{k}$ through model from (\ref{eqn:cableselection}) to (\ref{eqn:cableselection4})}}
\STATE{\textit{Formulate MILP model from (\ref{eqn:objective}) to (\ref{eqn:c8})}}
\IF{\textit{($k_{o}=1$)}}
    \STATE{\textit{Warm start with initial feasible point $\boldsymbol{I}$}}
\ELSE
    \STATE{\textit{Warm start with feasible point obtained in $\boldsymbol{O}_{k_{o}-1}$}}
\ENDIF
\STATE{\textit{Solve MILP model from (\ref{eqn:objective}) to (\ref{eqn:c8})}}
\STATE{\textit{$\boldsymbol{Z}_{k_{o}} = \left\lbrace(i,j)\right\rbrace:x_{ij}=1 \quad \forall (i,j) \in \boldsymbol{A^{'}_{r}}$}}
    \IF{\textit{($k_{o}>1$)}}
        \IF{\textit{($\boldsymbol{Z}_{k_{o}} \subset \boldsymbol{\Gamma}_{k_{o}-1}$)}}
           \STATE{\textit{Save best feasible point found: $\boldsymbol{O}_{k_{o}} = x_{ij} \cup y^k_{ij} \quad \forall (i,j) \in \boldsymbol{A^{'}_{r}} \wedge k \in \left\lbrace1,\cdots,f(i)\right\rbrace$}}
           \STATE{\textit{Break}}
        \ELSE
            \STATE{\textit{Save feasible point found: $\boldsymbol{O}_{k_{o}} = x_{ij} \cup y^k_{ij} \quad \forall (i,j) \in \boldsymbol{A^{'}_{r}} \wedge k \in \left\lbrace1,\cdots,f(i)\right\rbrace$}}
            \STATE{\textit{$k_{o}\leftarrow k_{o}+1$}}
        \ENDIF
    \ENDIF
\ENDFOR
\STATE{\textit{Recalculate gaps}}
\end{algorithmic}
\rule[1ex]{0.95\linewidth}{0.5pt}
\end{algorithm}
From line 1 to 12 the task is to efficiently solve a  feasibility problem. The idea is to subsequently increase $\upsilon$ from an initial value $\upsilon=\upsilon_{f}=\upsilon_{f_{min}}$ to a maximum value $\upsilon=\upsilon_{f}=\upsilon_{f_{max}}$, with steps $\upsilon_{f_{\delta}}$, until a feasible point is found. If this is achieved in iteration $k_f$, the first task is terminated with a feasible point $\boldsymbol{I}$. Conversely, if the model is infeasible, the candidate arcs set is augmented with $\upsilon_{f_{\delta}}$ units, and the process is taken to the iteration $k_f+1$, where a new trial is attempted. \\In order to formulate the MILP model, the cost coefficients calculation from (\ref{eqn:cableselection}) to (\ref{eqn:cableselection4}) is omitted by setting them equal to zero, and the black-box MILP solver terminates when the first feasible point is found. \\At this point, the Algorithm \ref{alg:al2} requires as parameters $\upsilon_{f_{min}}$, $\upsilon_{f_{\delta}}$, and $\upsilon_{f_{max}}$. The greater $\upsilon_{f_{min}}$ and $\upsilon_{f_{\delta}}$ the less efficient the feasibility problem becomes, however, increasing the odds to defining a feasible instance of the problem promptly. \\Likewise, from line 13 to 35, the global optimization task is performed. The target is to obtain a feasible point with a given relative gap $\epsilon$, expressed as the relative difference of the best feasible point ($\tau$) minus the best achievable value objective ($\kappa$), with respect to $\tau$. These values are indexed by iteration number. Similarly to the feasibility task, the iterative process increases the candidate arcs set from $\upsilon=\upsilon_{o}=\upsilon_{o_{min}}$ to $\upsilon=\upsilon_{o}=\upsilon_{o_{max}}$, with steps $\upsilon_{o_{\delta}}$. \\The termination criterion is when the set $\boldsymbol{Z}_{k_{o}}$ of active variables $x_{ij}=1$ of the problem defined in the iteration $k_o$, is a subset of the arcs set $\boldsymbol{A^{'}_{r}}$ defined in the previous iteration $k_{o}-1$ ($\boldsymbol{\Gamma}_{k_{o}-1}$). \\In this way, it is inferred that it is not longer necessary to increase $\upsilon_{o}$, as the optimum variables have been already provided in the previous iteration. \\To guarantee along the process a monotonously decreasing value of the objective function, in iteration $k_o$, the mathematical model is warm-started with the feasible solution found in $k_{o}-1$ ($\boldsymbol{O}_{k_{o}-1}$). This strategy may help in shortening the convergence time for the sub-instance $k_o$. \\Conceptually, Algorithm \ref{alg:al2} intends to determine a reduced search space, where the global minimum point is hopefully included. If only one reduced problem was solved given a $\upsilon$, it would not be possible  to infer about the quality of the solution, and the calculated gap for that particular instance could not represent the global domain of the full problem, potentially leading to an overestimation.\\ For the global optimization task, Algorithm \ref{alg:al2} requires as parameters $\upsilon_{o_{min}}$, $\upsilon_{o_{\delta}}$, and $\upsilon_{o_{max}}$ for the global optimization task. Naturally, $\upsilon_{o_{min}} \geq \upsilon_{f_{max}}$, and it is reasonable to consider $\upsilon_{o_{\delta}}>\upsilon_{f_{\delta}}$. By proper adjustment of the previous parameter, in best case scenario, the full Algorithm is concluded for $k_f=1$ and $k_o=2$. \\Although for every iteration the maximum required gap $\epsilon$ is equally fixed, the equivalent calculated gap, having as reference the full-size domain, varies. Larger values of $\upsilon_{k_{o}}$ lead to equal or lower values of $\kappa_{k_{o}}$. This causes that in general, $\tau_{k_{o}}$ is also lower, until the ideal reduced search space is found, when equal values of $\tau_{k_{o}}$ should be obtained. \\Therefore, after the termination of the algorithm, a gap updating procedure is performed based on the last calculated value of $\kappa_{k_{o}}$, to recalculate the relative difference for all previous iterations respect to this value (line 36). Let the recalculated gap in the global iteration $k_g$, including the feasibility and global optimization problems, be $\epsilon_{k_g}$. In this sense, an evolution of the gap in function of the iterations is available, providing further insights and the sense of convergence, as the objective value decreases monotonically.
\section{Computational Experiments}
\label{experiments}
The following experiments have been carried out on an Intel Core i7-6600U CPU running at 2.50 GHz and with 16 GB of RAM. The chosen MILP solver is the branch-and-cut solver implemented in IBM ILOG CPLEX Optimization Studio V12.7.1 \cite{IBM2015IBMManual}. \\In the Section \ref{sensitivity} a sensitivity analysis of the parameters $\upsilon_{o_{min}}$, and $\upsilon_{o_{\delta}}$ is presented. This is achieved by applying the proposed methodology of Section \ref{framework}, on two OWFs using several sets of parameters. \\In Section \ref{benchmarking}, the proposed method is benchmarked against the results obtained through a different approach published in the scientific literature \cite{Fischetti2018OptimizingLosses}. With this aim, the same testbed is employed, while assuming the same considerations, such as, objective (\ref{eqn:cableselection})-(\ref{eqn:cableselection4}), and constraints embodied by the equations (\ref{eqn:c0})-(\ref{eqn:c8}).
\subsection{Algorithm's Parameters Sensitivity Analysis}\label{sensitivity}
The two real OWFs West of Duddon Sands (WDS) and Thanet (TH) are used for the parameter sensitivity analysis. The information regarding these OWFs is provided in \cite{ESCAKIS-ORCA}. \\Table \ref{tab:param_inputs} displays the main parameters for the sensitivity analysis. In the case of WDS, the objective function is defined by a combined total economic cost including the initial investment, and the total electrical power losses of the cable layout (IP). Whilst for TH, the target is to minimize the initial investment of the cable layout, defined as the cables capital expenditures, and installation costs. \\The difference between objective functions is not relevant for the purposes of the sensitivity analysis, but rather, to show the capability of the method to support both cases. Economic parameters for the discounted cash flows calculations are also presented.\\
Both instances are challenging to solve given their large size (with more than $100$ WTs each), the large number of cables sizes considered (up to three cables types), and maximum capacity $U=\max{\boldsymbol{U}}$, spanning from $10$ to $13$ (see Table \ref{tab:param_inputs}). High-level information such as WT power and number ($P_n$ and $n_w$), OSSs number ($n_o$), and the maximum number of allowed feeders connected to the OSSs ($\phi$) is available. The limit $\phi$ provides a hard binding constraint. It is fixed to a practical value usually considered by OWFs developers.\\
Finally, the rest of parameters for the optimization are displayed. The collection systems voltage level ($V_n$) is the traditionally used, and two different set of cables are considered. The electrical and thermal information of the cables is available in \cite{ABB2018XLPEGuide}. \\The capacity set $\boldsymbol{U}$ is determined for each instance, similarly, for the total cost per kilometer, including capital expenditures ($\boldsymbol{C_{c}}$), and installation costs $\boldsymbol{C_{p}}$. The cables $\left\lbrace1,2,3\right\rbrace$ have cross-sections of {\SI{240}{mm^2}, \SI{500}{mm^2}, and \SI{1000}{mm^2}}, respectively. \\The set represented by $\left\lbrace4,5\right\rbrace$ is used for TH, and it has associated electrical and economic parameters matching those considered in the benchmark work \cite{Fischetti2018OptimizingLosses}, while neglecting the capacitive currents (see Section \ref{benchmarking} for benchmark analysis).\\ 
The whole framework compacted in Algorithm \ref{alg:al2} intends to find an approximation of the minimum search space of a given instance. It is presumed the finding the model best feasible point (or near to it) for a given required maximum gap, having as reference the full (global) problem size. \\Each problem instance is composed by $k_g=k_f+k_o$ sub-instances (or general iterations), which in turn include those iterations required for solving the feasibility ($k_f$), and global optimization serial problems ($k_o$); this implies that, in the best-case scenario, a maximum number of iterations equal to $k_f=1$ and $k_o=2$ should be enough to solve any problem instance.\\ However, factors such as the adequate setting of the algorithm's parameters and the specific spatial characteristics influence the number of global iterations $k_g$, and the circumvention of local minima. \\The avoidance of local minima is crucial, as it has a lifetime impact over a project. Successful implementations depend on the chosen values of $\upsilon_{o_{min}}$, and $\upsilon_{o_{\delta}}$.
\begin{table*} [!t]
\centering
\captionsetup{justification=centering, labelsep=newline,textfont = sc}
\caption{Main inputs parameters for sensitivity analysis}
\label{tab:param_inputs}
\setlength{\tabcolsep}{3pt}
\begin{tabular}{S[table-number-alignment = center,table-column-width=0.50cm] S[table-number-alignment = center,table-column-width=0.50cm] S[table-number-alignment = center,table-column-width=0.73cm] S[table-number-alignment = center,table-column-width=1.53cm] S[table-number-alignment = center,table-column-width=0.53cm] : S[table-number-alignment = center,table-column-width=1.53cm] S[table-number-alignment = center,table-column-width=0.52cm] S[table-number-alignment = center,table-column-width=0.72cm] S[table-number-alignment = center,table-column-width=0.72cm] S[table-number-alignment = center,table-column-width=1.42cm] : S[table-number-alignment = center,table-column-width=1.42cm] S[table-number-alignment = center,table-column-width=1.32cm] S[table-number-alignment = center,table-column-width=1.32cm] S[table-number-alignment = center,table-column-width=2.45cm]} 
\toprule
{OWF} & {Obj.} & $r$ [$\%$] & ${c}_{e}$ [\EUR{}/\si{MWh}] & {$m$} & {$P_n$ [\si{MW}]} & {$n_w$} & {$n_o$} & {$\phi$} & {$\eta$} & {$V_n$ [{\si{kV}}]} & {$\boldsymbol{T}$} & {$\boldsymbol{U}$} & {$\boldsymbol{C_{c}}+\boldsymbol{C_{p}}$ [M\EUR{}/\si{\km}]}\\ 
\midrule
{WDS} & {IP} & 5 & 40 & 30 & 3.6 & 108 & 1 & 10 & 1 & 33 & {$\left\lbrace1,2,3\right\rbrace$} & {$\left\lbrace7,10,13\right\rbrace$} & {$\left\lbrace0.36,0.58,0.90\right\rbrace$}\\
{TH}  & {I}  & {-} & {-}  & {-}  & 3 & 100 & 1 & 10  & 1 & 33 & {$\left\lbrace4,5\right\rbrace$} & {$\left\lbrace7,10\right\rbrace$} & {$\left\lbrace0.44,0.62\right\rbrace$} \\
\bottomrule
\end{tabular}
\end{table*}
\begin{figure*}[!t]
\centering
  \subfloat[Sensitivity Analysis for West of Duddon Sands (WDS).]{\includegraphics[width=0.4\textwidth]{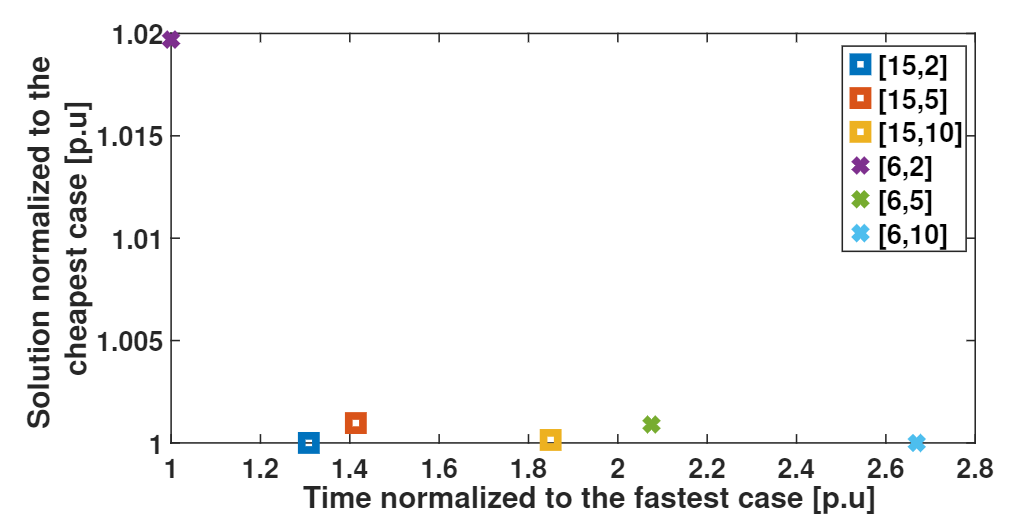}\label{fig:senwds}}
  \hfill
  \subfloat[Sensitivity Analysis for Thanet (TH).]{\includegraphics[width=0.4\textwidth]{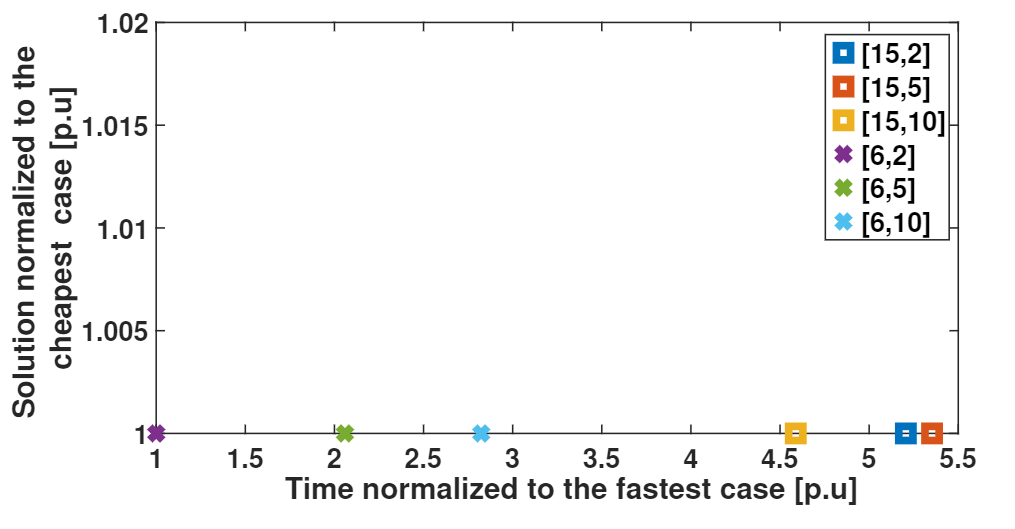}\label{fig:senth}}
\caption{Sensitivity Analysis Results.}
\label{fig:senwdsth}
\end{figure*}
\begin{figure*}[!t]
\centering
  \subfloat[West of Duddon Sands (WDS) Designed Cable Layout.  Yellow line: Cable \SI{1000}{mm^2}, Orange line: Cable \SI{500}{mm^2}, Blue line: \SI{240}{mm^2}.]{\includegraphics[width=0.4\textwidth]{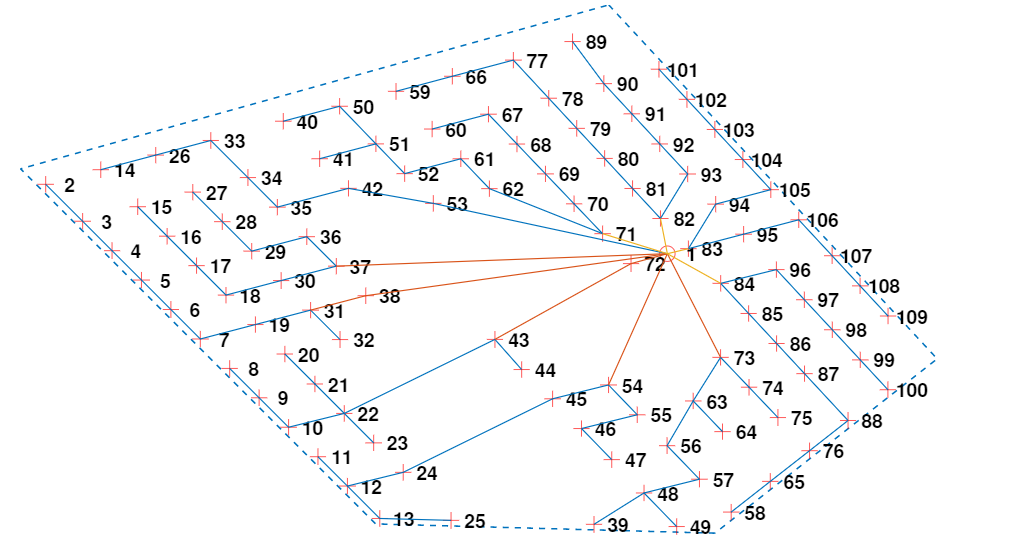}\label{fig:wds}}
  \hfill
  \subfloat[Thanet (TH) Designed Cable Layout.  Orange line: Cable supporting up to 10 WTs, Blue line: Cable supporting up to 7 WTs.]{\includegraphics[width=0.4\textwidth]{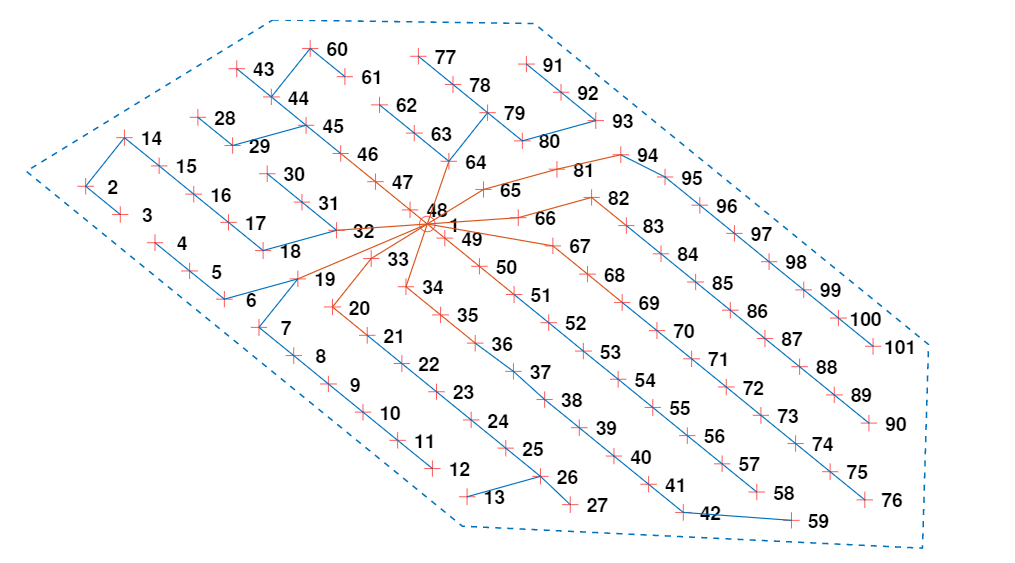}\label{fig:TH}}
\caption{Designed Cable Layouts,}
\label{fig:wdsandla}
\end{figure*}
The sensitivity analysis elaborates on how the iterative reduction of the search space through Algorithm \ref{alg:al2} can lead to sub-optimal solutions. \\In pursuance of this, two values associated to the initial value parameter, $\upsilon_{o_{min}}$, are examined: $\{6,15\}$; while the increase steps parameter, $\upsilon_{o_{\delta}}$, is set between $\{2,5,10\}$. All possible combinations of both parameters are applied to WDS and TH OWFs through the Algorithm \ref{alg:al2}.\\
The results are graphically shown in Figure \ref{fig:senwdsth} presenting quality of the solution (ordinate) versus computing time (abscissa). Both are normalized to the best corresponding case. Each case is depicted in the legend with [$\upsilon_{o_{min}}$,$\upsilon_{o_{\delta}}$], i.e. [6,2] and so on.\\The results show that for the WDS OWF, the solution for $\upsilon_{o_{min}}=6$, and $\upsilon_{o_{\delta}}=2$ is nearly $2\%$ more expensive than the other cases (see Figure \ref{fig:senwds}, case [6,2]). For the rest of parameters sets, the obtained solutions are, in practical terms, the same. The costs differences among them are due to the inclusion of total electrical power losses, which are sensitive to small changes of cables lengths. The computing times varies almost in a positively correlated fashion with the values of these parameters, with the fastest case being also the most expensive.\\ 
This could be explained by the particular spatial characteristics of WDS, which impact the number of minimum candidate arcs to cover the global minimum. The optimized layout shown in Figure \ref{fig:wds} ($\upsilon_{o_{min}}=15$ and $\upsilon_{o_{\delta}}=5$) evidentiate the non-uniform distribution of WTs in the plane, with empty areas around WTs number 43 and number 44, for instance. This can be interpreted as a higher degree of freedom in the design, having the possibility to interconnect WTs located further away from each other. \\One should note that the optimized layout has a connection between WTs number 72 and 43, the latter being out of the top-15 closest WTs for the former (in fact, is the 17th closest), hence, making it necessary to increase $\upsilon_o$ to $20$, to include this connection.\\
Returning to the case of $\upsilon_{o_{min}}=6$, and $\upsilon_{o_{\delta}}=2$, the step increase value is not enough to provide a significantly broader search space to improve the solution quality. \\The results suggest that the values of $\upsilon_{o_{\delta}}=5$, and $\upsilon_{o_{\delta}}=10$, result in an improvement of the solution in each iteration, and confirm in the last iteration the covering of minimum candidate arcs set. A value of $\upsilon_{o_{\delta}}=2$ in this case provides the same output, but evidently with potentially higher chances of falling into suboptimal points for other OWFs. \\
For the TH OWF, all the cases generate the same result, as shown in Figure \ref{fig:senth}. In contrast to WDS OWF, Thanet presents a regular (grid-based) micrositing layout  (see Figure \ref{fig:TH}), hence, less degree of flexibility, which translates in an optimum point represented by WTs each connected to maximum the 6th closest WT. \\In order to increase the likelihood for getting a balance between solution quality, computing time, and memory requirements, a heuristic rule of considering $\upsilon_{o_{min}}=15$, and $\upsilon_{o_{\delta}}=5$ is implemented in this manuscript.  \\Certainly, every single OWF should be individually analysed to come up with tailored parameters, but results point out that these settings may lead to adequate terminations in most of the problems.\\Table \ref{tab:det_results} summarizes in-detail the results of the experiments using the proposed framework with $\upsilon_{f_{min}}=5$, $\upsilon_{f_{\delta}}=1$, $\upsilon_{f_{max}}=15$, $\upsilon_{o_{min}}=15$, $\upsilon_{o_{\delta}}=5$, and $\upsilon_{o_{max}}=50$. \\For each of the sub-instances, the output CT1 includes the total processing time in order to generate the Constraints (\ref{eqn:c0}) to (\ref{eqn:c8}), and additionally, the solution of the independent sub-problems as defined in the model from (\ref{eqn:cableselection}) to (\ref{eqn:cableselection4}). Note that in the case of the feasibility problems the latest procedure is circumvented by fixing the cost coefficients equal to zero. \\Likewise, CT2 is the computing time to solve the main mathematical model (from (\ref{eqn:objective}) to (\ref{eqn:c8})) for a given maximum gap, in the case of WDS $\epsilon=0.5\%$, and for TH $\epsilon=0.3\%$. An additional experiment with $\epsilon=0.5\%$ for TH has been run, with a duration of around an hour, finding a solution only $0.07\%$ more expensive than the one presented in Table \ref{tab:det_results}, with a calculation time approx. $50\%$ smaller. In principle, any imposed maximum gap is supported, at the potential expense of a very steep increase on computing time.
\begin{table} [!t]
\centering
\captionsetup{justification=centering, labelsep=newline,textfont = sc}
\caption{In-detail results summary}
\label{tab:det_results}
\begin{tabular}{S[table-number-alignment = center,table-column-width=0.1cm] S[table-number-alignment = center,table-column-width=0.13cm] S[table-number-alignment = center,table-column-width=0.13cm] S[table-number-alignment = center,table-column-width=0.13cm] S[table-number-alignment = center,table-column-width=0.13cm] S[table-number-alignment = center,table-column-width=0.32cm] S[table-number-alignment = center,table-column-width=0.32cm] S[table-number-alignment = center,table-column-width=0.32cm] S[table-number-alignment = center,table-column-width=0.32cm] S[table-number-alignment = center,table-column-width=0.63cm] S[table-number-alignment = center,table-column-width=0.63cm] S[table-number-alignment = center,table-column-width=0.63cm]}
\toprule
 & \multicolumn{2}{c}{{$k_f$/}} & \multicolumn{2}{c}{{$\upsilon_{f}$/}} & \multicolumn{2}{c}{{CT1}} & \multicolumn{2}{c}{{CT2}} & {Obj} & {$\epsilon_{k_{f}}$/$\epsilon_{k_{o}}$} & {$\epsilon_{k_g}$} \\ 
 & \multicolumn{2}{c}{{$k_o$}} &\multicolumn{2}{c}{{$\upsilon_{o}$}} &{Time}&{Unit} & {Time}&{Unit} & {[M\EUR{}]} & {[$\%$]} & {[$\%$]} \\
\midrule
           {W}            & {$k_f\leftarrow$} & 1 & {$v_f\leftarrow$} & 5 & 7.89 & {s} & 36.93 & {s} & 72.84 & 0 & 38.62 \\
           {D}           & {$k_o\leftarrow$} & 1 & {$v_o\leftarrow$} & 15 & 0.08 & {h} & 1.01 & {h} & 45.70 & 0.46 & 2.17 \\
           {S}               & {$k_o\leftarrow$} & 2 & {$v_o\leftarrow$} & 20 & 0.11 & {h} & 0.36 & {h} & 44.93 & 0.47 & 0.50 \\
                       & {$k_o\leftarrow$} & 3 & {$v_o\leftarrow$} & 25 & 0.14 & {h} & 0.16 & {h} & 44.93 & 0.50 & 0.50 \\
\midrule                     
{T} & {$k_f\leftarrow$} & 1 & {$v_f\leftarrow$} & 5 & 6.65 & {s} & 10.37 & {s} & 40 & 0 & 33.60 \\
{H}                      & {$k_o\leftarrow$} & 1 & {$v_o\leftarrow$} & 15 & 15.80 & {s} & 1.02 & {h} & 26.64 & 0.30 & 0.30 \\
                      & {$k_o\leftarrow$} & 2 & {$v_o\leftarrow$} & 20 & 30.22 & {s} & 0.75 & {h} & 26.64 & 0.30 & 0.30 \\     
\bottomrule
\end{tabular}
\end{table}
The column Obj presents the objective value of the best feasible point obtained under the explained conditions, after the termination of each sub-instance calculation. Finally, the columns $\epsilon_{k_{f}}$/$\epsilon_{k_{o}}$, and $\epsilon_{k_g}$ gather the calculated gap of the best feasible point, and the recalculated gap for each sub-instances, respectively, after the Algorithm \ref{alg:al2} terminates as explained in Section \ref{thealgorithm}. \\By means of this strategy, is surmised the delimitation of a reduced search space representative of the global problem, including the global minimum point. \\In the particular case of WDS, the results say that an unnecessary continuation to a fifth cycle, i.e. $k_g=5$, was avoided due to the warm start strategy. The objective value does not change between $k_o=2$ and $k_o=3$, this being a clear (but not mathematically conclusive) signal of the successful convergence of the method. \\Furthermore, an overestimation of $\epsilon_{k_{o}}$ is seen for $v_{o}=15$ when referenced to the approximated minimum search space, as $\epsilon_{k_g}=2.17\%>0.5\%$. It should be noted that the first feasible point is found in almost $45$ seconds, and from there, the best feasible point is found in very reasonable computing time. A reduction of the gap from $38.62\%$ to $0.50\%$ is obtained in only $1$ h:$53$ min (one hour and $53$ minutes).\\
In contrast to WDS, only three global iterations are required for TH, mainly due to the more uniform distribution of WTs, as shown in the Figure \ref{fig:TH}. For TH OWF, a feasible point is obtained in $17$ s, and in $1$ h:$47$ min the gap is improved from $33.60\%$ to $0.30\%$.\\The proposed procedure can be impacted by multiple valid solutions around the required maximum gap. Nevertheless, this would cause only false continuations of the Algorithm \ref{alg:al2}, rather than affect the solution quality. The previous point is partially palliated by the likely decrease of computing time in subsequent iterations, given the warm start point provided from the previous step. If the objective value was used as criterion to stop the algorithm, a maximum ratio of objective change per subsequent iterations should be set, which would be open to different assessment criteria, and potentially could lead to false terminations.
\\For both WDS and TH, the values $\epsilon_{k_g}$ reported in Table \ref{tab:det_results} for the last iteration, are representative of the full problem; fact validated when solving them without reduced search space (full arcs set), provided the best feasible point available.
\subsection{Benchmarking}
\label{benchmarking}
A testbed of real-world cases, presented in Table \ref{tab:param_bench}, is employed. The OWFs names are given in the Acronyms definition. Results comparison are available in Table \ref{tab:res_bench}.\\
\begin{table*} [!t]
\centering
\captionsetup{justification=centering, labelsep=newline,textfont = sc}
\caption{Main inputs parameters for benchmarking}
\label{tab:param_bench}
\setlength{\tabcolsep}{3pt}
\begin{tabular}{S[table-number-alignment = center,table-column-width=0.7cm] S[table-number-alignment = center,table-column-width=0.7cm] S[table-number-alignment = center,table-column-width=0.7cm] S[table-number-alignment = center,table-column-width=0.7cm] S[table-number-alignment = center,table-column-width=1.5cm] S[table-number-alignment = center,table-column-width=0.5cm] : S[table-number-alignment = center,table-column-width=1.2cm] S[table-number-alignment = center,table-column-width=0.7cm] S[table-number-alignment = center,table-column-width=0.7cm] S[table-number-alignment = center,table-column-width=0.7cm] S[table-number-alignment = center,table-column-width=1cm] : S[table-number-alignment = center,table-column-width=1cm] S[table-number-alignment = center,table-column-width=1.3cm] S[table-number-alignment = center,table-column-width=1.3cm] S[table-number-alignment = center,table-column-width=2.2cm]} 
\toprule
{Ins.} & {OWF} & {Obj.} & $r$ [$\%$] & ${c}_{e}$ [\EUR{}/\si{MWh}] & {$m$} & {$P_n$ [\si{MW}]} & {$n_w$} & {$n_o$} & {$\phi$} & {$\eta$} & {$V_n$ [{\si{kV}}]} & {$\boldsymbol{T}$} & {$\boldsymbol{U}$} & {$\boldsymbol{C_{c}}+\boldsymbol{C_{p}}$ [M\EUR{}/\si{\km}]}\\ 
\midrule
{1} & {HR1}  & {I} & {-} & {-} & {-} & 2 & 80 & 1 & 10 & 1 & 33 & {$\left\lbrace6,7,8\right\rbrace$} & {$\left\lbrace7,11,13\right\rbrace$} & {$\left\lbrace0.37,0.39,0.43\right\rbrace$}\\
{2} & {HR1}  & {I} & {-} & {-} & {-} & 2 & 80 & 1 & 10 & 1 & 33 & {$\left\lbrace9,10\right\rbrace$} & {$\left\lbrace7,12\right\rbrace$} & {$\left\lbrace0.44,0.45\right\rbrace$}\\
{3} & {HR1}  & {I} & {-} & {-} & {-} & 2 & 80 & 1 & 10 & 1 & 33 & {$\left\lbrace4,5\right\rbrace$} & {$\left\lbrace10,14\right\rbrace$} & {$\left\lbrace0.44,0.62\right\rbrace$}\\
{4} & {O}  & {I} & {-} & {-} & {-} & 5 & 30 & 1 & 4 & 1 & 33 & {$\left\lbrace11,12\right\rbrace$} & {$\left\lbrace5,10\right\rbrace$} & {$\left\lbrace0.41,0.61\right\rbrace$}\\
{5} & {O}  & {I} & {-} & {-} & {-} & 5 & 30 & 1 & 4 & 1 & 33 & {$\left\lbrace13,14\right\rbrace$} & {$\left\lbrace4,9\right\rbrace$} & {$\left\lbrace0.38,0.63\right\rbrace$}\\
{6} & {DT}  & {I} & {-} & {-} & {-} & 3.6 & 80 & 1 & 10 & 1 & 33 & {$\left\lbrace6,7,8\right\rbrace$} & {$\left\lbrace4,6,8\right\rbrace$} & {$\left\lbrace0.37,0.39,0.43\right\rbrace$}\\
{7} & {DT}  & {I} & {-} & {-} & {-} & 3.6 & 80 & 1 & 10 & 1 & 33 & {$\left\lbrace4,5\right\rbrace$} & {$\left\lbrace6,8\right\rbrace$} & {$\left\lbrace0.44,0.62\right\rbrace$}\\
{8} & {TH}  & {I}  & {-} & {-}  & {-}  & 3 & 100 & 1 & 10  & 1 & 33 & {$\left\lbrace13,14\right\rbrace$} & {$\left\lbrace7,15\right\rbrace$} & {$\left\lbrace0.38,0.63\right\rbrace$} \\
{9} & {TH}  & {I}  & {-} & {-}  & {-}  & 3 & 100 & 1 & 10  & 1 & 33 & {$\left\lbrace4,5\right\rbrace$} & {$\left\lbrace7,10\right\rbrace$} & {$\left\lbrace0.44,0.62\right\rbrace$} \\
\cdashline{1-15}
{10} & {LA}  & {IP} & 5 & 40 & 30 & 3.6 & 175 & 2 & 10 & 1 & 33 & {$\left\lbrace1,2,3\right\rbrace$} & {$\left\lbrace7,10,13\right\rbrace$} & {$\left\lbrace0.36,0.58,0.90\right\rbrace$}\\
\bottomrule
\end{tabular}
\end{table*}
\begin{table*} [!t]
\centering
\captionsetup{justification=centering, labelsep=newline,textfont = sc}
\caption{Results for benchmarking}
\label{tab:res_bench}
\setlength{\tabcolsep}{3pt}
\begin{tabular}{S[table-number-alignment = center,table-column-width=0.95cm] S[table-number-alignment = center,table-column-width=0.95cm] S[table-number-alignment = center,table-column-width=1.5cm] S[table-number-alignment = center,table-column-width=1.5cm] S[table-number-alignment = center,table-column-width=1.5cm] S[table-number-alignment = center,table-column-width=1.5cm] S[table-number-alignment = center,table-column-width=1.5cm] S[table-number-alignment = center,table-column-width=1.5cm] S[table-number-alignment = center,table-column-width=1.5cm] S[table-number-alignment = center,table-column-width=1.5cm] S[table-number-alignment = center,table-column-width=1.8cm]} 
\toprule
{Ins.} & {OWF} & {Obj.\cite{Fischetti2018OptimizingLosses} [M\EUR{}]} & {Obj. [M\EUR{}]} & {Diff.[M\EUR{}]} & {{$\epsilon_{k_g}$}\cite{Fischetti2018OptimizingLosses} [\%]} & {{$\epsilon_{k_g}$} [\%]} & {Diff.[\%]}  & {Time\cite{Fischetti2018OptimizingLosses} [min}] & {Time [min}] & {Diff.[min]}\\
\midrule
{1} & {HR1} & 19.44 & 19.44 & 0 & 0.01 & 0.01 & 0 & 30 & 1.57 & 28.43 \\
{2} & {HR1} & 22.61 & 22.61 & 0 & 0.01 & 0.01 & 0 & 30 & 1 & 29 \\
{3} & {HR1} & 23.48 & 23.48 & 0 & 0.17 & 0.01 & 0.16 & 1440 & 7.26 & 1432.74 \\
{4} & {O} & 8.05 & 8.05& 0 & 0 & 0.01 & -0.01 & 0.42 & 0.46 & -0.04 \\
{5} & {O} & 8.36 & 8.36& 0 & 0 & 0.01 & -0.01 & 1.95 & 0.45 & 1.50 \\
{6} & {DT} & 38.98 & 38.98 & 0 & 4.85 & 0.01 & 4.84 & 30 & 1.88 & 28.12 \\
{7} & {DT} & 50.38 & 49.83 & 0.55 & 7.36 & 0.01 & 7.35 & 10 & 2.65 & 7.35 \\
{8} & {TH} & 22.34 & 22.31 & 0.03 & 3.37 & 0.7 & 2.67 & 5 & 627.60 & -622.60 \\
{9} & {TH} & 26.64 & 26.64 & 0 & 2.51 & 0.3 & 2.21 & 1440 & 107.02 & 1332.98\\
\cdashline{1-11}
{10} & {LA} & {-} & 68.97 & {-} & {-} & 0.5 & {-} & {-} & 748 & {-} \\
\bottomrule
\end{tabular}
\end{table*}
This testbed has been mostly extracted from \cite{Fischetti2018OptimizingLosses}, by selecting the most challenging instances. In \cite{Fischetti2018OptimizingLosses} a significantly different approach (different mathematical formulation, for instance) from the one proposed in this article, consisting in a matheuristic model, has been designed, implemented, and tested. \\The results in \cite{Fischetti2018OptimizingLosses} where obtained with computational resources similar to those used in this manuscript (IX CPU X5550 running at 2.67GHz, CPLEX 12.6). \\In order to provide a fair comparison, all practical and technical constraints are conceptually equivalent, while the objective function is the initial investment, as losses are computed differently. Likewise, capacitive currents have been neglected as they are not included in \cite{Fischetti2018OptimizingLosses}.
\\The comparison results are presented in Table \ref{tab:res_bench} (instances 1-9). A small (O), two large (HR1 and DT), and a very large OWF (TH) are studied. Each instance is defined by a OWF, and a set of cables available $\boldsymbol{T}$. \\Three aspects are compared: (i) solution quality, (ii) calculated gap, and (iii) computing time; each of them are directly compared by inspecting the columns Diff.[M\EUR{}], Diff.[\%], and Diff.[min], respectively. In all metrics, positive value means better performance for the method proposed in this article. 
\\Regarding solution quality, it can be seen that, for all instances, the obtained solutions are equal or lower than in the benchmark work. For instance nr.7, around \EUR{550,000}, representing around $1\%$ of the total cost, are saved. Particularly for this instance, the gap is improved from $7.36\%$ to $0.01\%$, while  simultaneously reducing the computing time. \\When the calculated gap in both methods is lower than $0.01\%$, the objective values are essentially the same. This validates that the dual values between the two models are also equivalent. \\In any instance the proposed method provides equally tight or even tighter solutions. The gap values reported in \cite{Fischetti2018OptimizingLosses} have been recalculated in this manuscript using the best feasible point as reference (instead of the best dual value obtained in $24$ h), to make them comparable with the proposed approach. \\Finally, in almost all the instances, the computing time is shorter, with exception of instance nr.8, where a considerable difference of $622$ min is observed. \\It is important to clarify that the reported computing times of the proposed method are for the whole running of Algorithm \ref{alg:al2}, including in all instances the final iteration $k_o=2$, necessary to confirm finding the global point; in contrast, in the benchmark work, the reported time is when the best feasible point has been found. Similarly, as mentioned before, $24$ h is the time limit to obtain the best dual value.\\For all the instances, the proposed method calculates feasible points in less than $40$ s, with recalculated gap of maximum $41\%$.
\\Besides the benchmark cases, an extra instance - not implemented in \cite{Fischetti2018OptimizingLosses} - is included to show the method's applicability to OWFs with multiple OSSs. LA OWF is the second largest project (measured by installed power) under operation, surpassed only by Walney Extension OWF (although with less number of WTs therefore potentially easier to solve). For this instance, an initial feasible point is obtained in only $3$ min, and the best feasible point is calculated in $12$ h:$28$ min, with the gap being improved by $46.87\%$.\\
\section{Conclusion}
\label{conclusion}
The proposed method provides a global optimization approach to solve the cable layout of OWFs collection systems. \\The main novelties of this manuscript are: (i) proposition of a model, able to provide very good solutions in very reasonable computing times, (ii) possibility to provide very tight solution quality certificates, (iii) integration of realistic and high-fidelity models to calculate total electrical power losses, and capacitive currents in the collection systems. \\An algorithmic framework for reducing the search space iteratively is the main technique used. The objective function supports the total economic costs, including initial investment, and lost revenues due to total electrical power losses in the project lifetime. \\The proposed methodology has been benchmarked against a state-of-the art method with significantly different approach (different MILP model, and application of up to four distinctive heuristics); with all practical and technical constraints conceptually equivalent. \\Ten real-world problem instances have been considered in the benchmark. The numerical results indicate that (i) the proposed algorithm provides, in general, at least equally good solutions, and in some cases, sizeable cheaper ones than the benchmark work, (ii) tighter gaps are calculated, in shorter computing times. \\The proposed algorithm also performs satisfactorily for large OWFs with multiple OSSs, where the clustering is intrinsically defined in the mathematical formulation. It does not require predecessor algorithms to group WTs into OSSs, avoiding in this way artificially biased solutions. 

\section*{Acknowledgment}
This research has received funding from the Baltic InteGrid Project (\url{www.baltic-integrid.eu/}). \\The authors thank Daniel Hermosilla Minguij\'on for his practical inputs throughout the development of the computational experiments, and Martina Fischetti, Lead Engineer - Operational Research in Vattenfall, for providing the data to perform the benchmark analysis.
\ifCLASSOPTIONcaptionsoff
  \newpage
\fi
\bibliographystyle{IEEEtran}
\bibliography{references}
\vskip -2\baselineskip plus -1fil
\begin{IEEEbiography}
[{\includegraphics[width=1in,height=1.25in,clip,keepaspectratio]{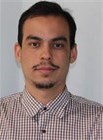}}]{Juan-Andr\'es P\'erez-R\'ua}  received the B.Sc.  degree in Electrical Engineering, with Summa Cum Laude distinction, from the Technological University of Bolivar, Colombia, in 2012, and the M.Sc. degree in Sustainable Transportation and Electrical Power Systems from the ISEC college, Coimbra, Portugal, the University of Nottinghan, England, and the University of Oviedo, Spain, in 2016. Currently, he is pursuing the Ph.D. degree in the Department of Wind Energy at the Technical University of Denmark (DTU). His present-day areas of interest include grid integration of renewables in power systems, and optimization.
\end{IEEEbiography}
\vskip -2\baselineskip plus -1fil
\begin{IEEEbiography}
[{\includegraphics[width=1in,height=1.25in,clip,keepaspectratio]{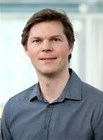}}]{Mathias Stolpe} received the Ph.D. degree in Optimization and Systems Theory from the Royal Institute of Technology (KTH), Stockholm, Sweden in 2003. He works as a Professor in the Department of Wind Energy at the Technical University of Denmark (DTU). His areas of research are structural and multidisciplinary optimization with focus on models and methods for global optimization.
\end{IEEEbiography}
\vskip -2\baselineskip plus -1fil
\begin{IEEEbiography}
[{\includegraphics[width=1in,height=1.25in,clip,keepaspectratio]{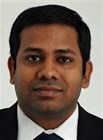}}]{Kaushik Das} (S'13, M'17) received the M.Tech. degree in Power system engineering from the IIT Kharagpur, India in 2011, and the Ph.D. degree from the Department of Wind Energy, Technical University of Denmark (DTU), Ris\o, Denmark in 2016. Currently, he is a Researcher in DTU Wind Energy. His research interest lies in hybrid power plants and grid integration of renewables in power systems. He is a member of IEA Wind, IEEE, CIGR\'{E} and other professional bodies.
\end{IEEEbiography}
\vskip -2\baselineskip plus -1fil
\begin{IEEEbiography}
[{\includegraphics[width=1in,height=1.25in,clip,keepaspectratio]{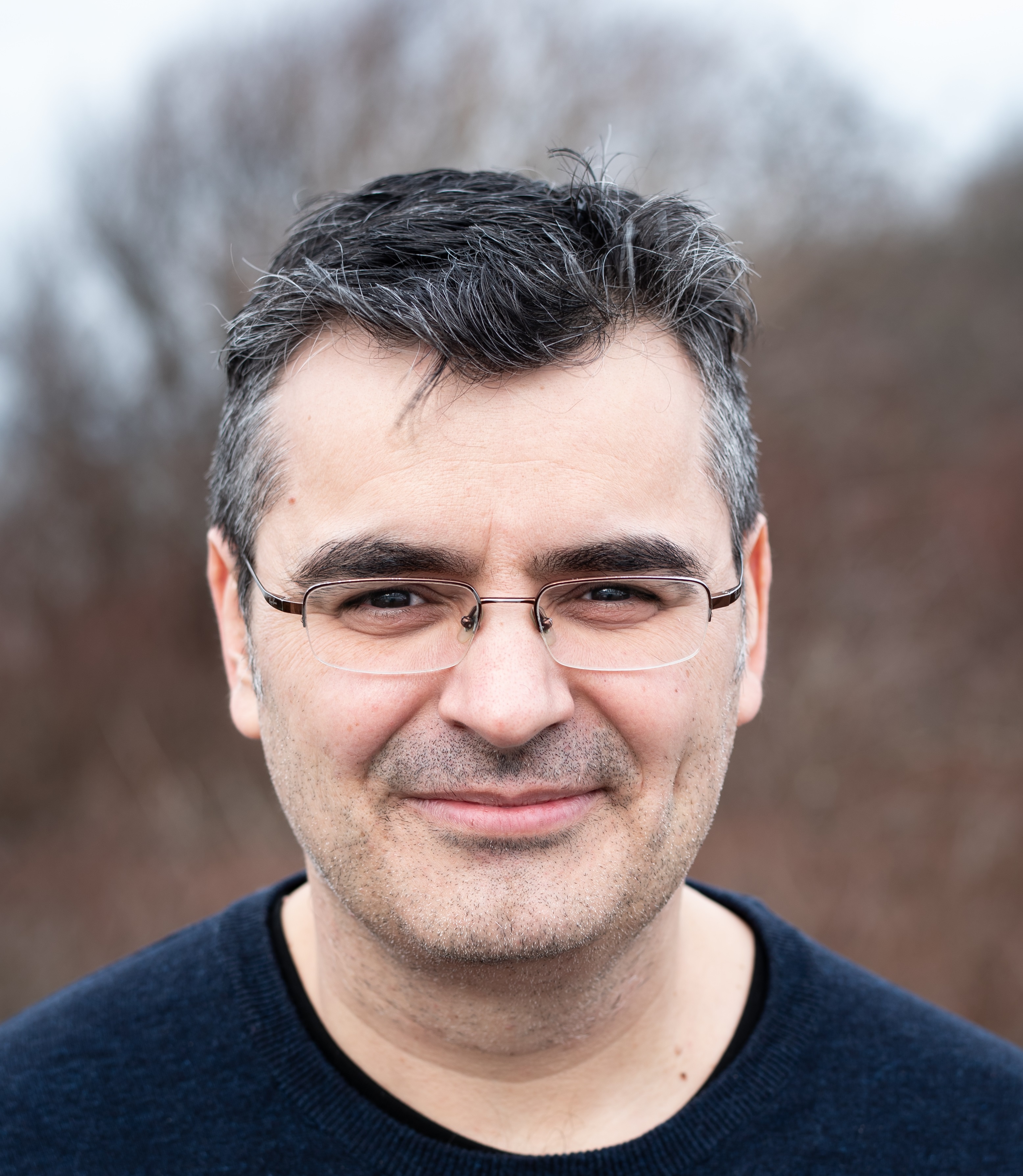}}]{Nicolaos A. Cutululis} (M'07, SM'18) received the M.Sc. and Ph.D. degrees, both in Automatic Control in 1998 and 2005, respectively. He is Professor in the Department of Wind Energy at the Technical University of Denmark. His main research interest is integration of wind power, with a special focus on offshore wind and grids. 
\end{IEEEbiography}
\end{document}